\documentclass[12pt,oneside]{article}
\usepackage{amsmath}
\usepackage{amsfonts}
\usepackage{amssymb}
\usepackage{graphicx}

\usepackage{latexsym}
\usepackage{amssymb, latexsym}
\usepackage{amsmath}
\usepackage{mathrsfs}
\usepackage{amsxtra}
\usepackage{amsthm}
\usepackage{setspace}
\usepackage{epsfig}
\usepackage{geometry}
\usepackage{graphicx}
\usepackage{blindtext}

\newtheorem{theorem}{Theorem}[section]
\newtheorem{corollary}[theorem]{Corollary}

\newtheorem{lemma}[theorem]{Lemma}

\newtheorem{claim}[theorem]{Claim}

\newtheorem{proposition}[theorem]{Proposition}

\theoremstyle{definition}
\newtheorem{definition}[theorem]{Definition}

\newtheorem{remark}[theorem]{Remark}

\date{}

\begin{document}

\author{Sergei Artamoshin
\thanks{The author wants to thank professor J\a'ozef Dodziuk for his interest to this paper and for many useful discussions.}\\ \\Department of Mathematical Science, Central Connecticut State University\\New Britain, USA\\ \\sartamoshin@gmail.com%
}

\title{\textbf{The Spherical Ratio of Two Points \\ and Its Integral Properties}}

\maketitle

\begin{abstract}
At the end of 19-th century, 1874, Hermann Schwartz found that for every point inside a planar disk, a two-dimensional Poisson Kernel can be written as a ratio of two segments, which he called as the geometric interpretation of that Kernel \cite{Schwarz}. About 90 years later, Lars Ahlfors, in his textbook, called this ratio as interesting \cite{Ahlfors}. We shall see here that every two different points, in a multidimensional space also define a similar ratio of two segments. The main goal of this paper is to study some integral properties of the ratio as well as introduce and proof One Radius Theorem for Spherical Ratio of two points. In particular, we introduce a new integration technique involving the ratio of two segments and a non-trivial integral equivalence leading to an integral relationship between Newtonian Potential and Poisson Kernel in any multi-dimensional space.
\end{abstract}

\bigskip

\section{Introduction}

\begin{figure}[!h]
    \centering
    \epsfig{figure=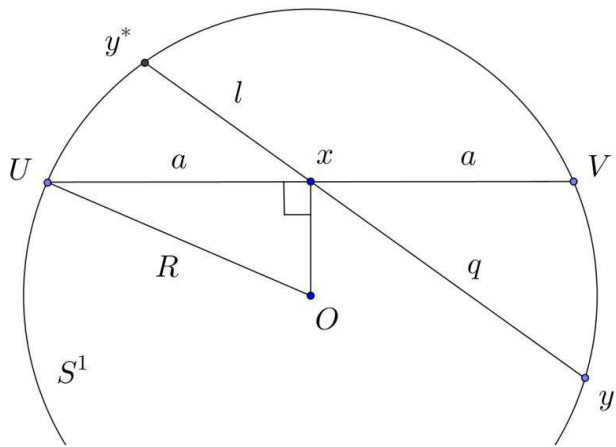,height=6cm}\\
    \caption{Hermann Schwarz Idea}\label{Geo_Inte_Schwarz}
\end{figure}

In this paper we introduce a new mathematical value $\omega$, called the spherical ratio of two points. Such ratio is a generalization of the ratio of two segments originally introduced by Hermann Schwartz at the end of 19-th century, (\cite{Schwarz}, pp. 359-361) or (\cite{Ahlfors}, pp. 168-169). While studying complex analysis, he found that for a point inside a planar disk, the Poisson Kernel can be written as a ratio of two segments. For the notations in the formula below, look at Figure~\ref{Geo_Inte_Schwarz}.

\begin{equation}\label{Chapter-goal}
    W(x,y)=\frac{R^2-|x|^2}{|x-y|^2}=\frac{a^2}{q^2}
    =\frac{lq}{q^2}=\frac{l}{q}\,.
\end{equation}

 We shall see that a similar ratio can be defined for every two different points in a multidimensional space and can serve as a useful tool for integration of some functions over sphere. The technique of integration together with the results obtained here can be applied for spectral analysis in a hyperbolic space. This happens because the ratio of two segments represents an eigenfunction of the hyperbolic Laplacian, while the averaging of $\omega$ over sphere represents a radial eigenfunction of the hyperbolic Laplacian. Therefore, all statements obtained in this paper can be restated in terms of radial eigenfunctions and used to obtain some results related, for example, to Dirichlet Eigenvalue Problem in a disc of constant negative curvature.

The main goal of this article is to describe the constraints for complex numbers $\alpha,\beta$ and for a fixed point $x\in\mathbb{R}^{k+1}\setminus\{S^k(R)\}$ under which the following equivalence holds.
\begin{equation}\label{Chapter-goal}
    \alpha+\beta=k\,\,\,\text{or}\,\,\,\alpha=\beta\quad\Longleftrightarrow\quad
    \int\limits_{S^k(R)}\omega^{\alpha}dS_y=
    \int\limits_{S^k(R)}\omega^{\beta}dS_y\,,
\end{equation}
where $\omega=\omega(x,y)$ is the spherical ratio of two points that is defined below.

\section{Definitions and Basic Results}

\begin{definition}[Spherical Ratio of Two Points]\label{Spherical_Ratio_Definition}
   Let $x,y\in\mathbb{R}^{k+1}$ and assume that $x\neq y$. Let $S^k(R)$ denotes the $k$-dimensional sphere of radius $R$ centered at the origin $O$. Then, let us introduce $x^*$ and $y^*$ as follows. If the line defined by $x$ and $y$ is tangent to $S^k(|x|)$, then we set $x^*=x$. Otherwise, $x^*$ be the point of $S^k(|x|)$ such that $x^*\neq x$ and $x^*, x, y$ are collinear. Note that $x^*=y$ in case $|x|=|y|$. Similarly, if the line through $x,y$ is tangent to $S^k(|y|)$, then we set $y^*=y$. Otherwise, $y^*$ be the point of $S^k(|y|)$ such that $y^*\neq y$ and $y^*, x, y$ are collinear. Note again that $y^*=x$ in the case $|x|=|y|$. Then, we can observe that $|x-y^*|=|x^*-y|$ and denote $$q=|x-y|\quad \text{and}\quad l=|x-y^*|=|x^*-y|\,.$$ We call the ratio of these two segments $l$ and $q$, as the spherical ratio of two points $x$ and $y$, i.e.,

     \begin{equation}\label{omega-definition}
        \omega(x,y)=\frac{l}{q}=\frac{|x-y^*|}{|x-y|}=\frac{|x^*-y|}{|x-y|}
        =\left|\frac{|x|^2-|y|^2}{|x-y|^2}\right|\,,
     \end{equation}
where the last expression is exactly the 2-D Poisson Kernel for a planar disk.
\end{definition}

\begin{corollary} Using the ratio of two segments $l$ and $q$ obtained above, we observe that $\omega(x,y)$ remains constant as long as $x$ stays at a circle tangent to $S^k(|y|)$ and centered at some point of the line $Oy$. This follows from the Lemma~\ref{omega_constant-lemma}, p.~\pageref{omega_constant-lemma} in Appendix. Then, the combination of Lemma~\ref{omega_constant-lemma} and Corollary~\ref{limit-for-omega} from p.~\pageref{limit-for-omega} shows that $\lim\limits_{x\rightarrow y}\omega(x,y)\in[0,\infty]$.
\end{corollary}

\begin{theorem}[\textbf{One Radius Theorem for} $\omega$]\label{Basic-theorem}

Let $S^k$,
$k\in\mathbb{N}$, be a $k$-dimensional sphere of radius $R$
centered at the origin $O$, $x, y\in\mathbb{R}^{k+1}$ such that $r=|x|\neq|y|=R$ and $\omega(x,y)$ is defined in \eqref{omega-definition}.
Until further notice we assume that $R$ and a point $x\notin S^k$ are fixed.
Then, the following Statements hold.

\begin{description}

  \item[(A)]\label{Sphe-pro-direct-Statement} If $\alpha, \beta\in\mathbb{C}$ and $\alpha+\beta=k$, then
\begin{equation}\label{Sphe-Pro-Statement-2}
    \int\limits_{S^k}\omega^\alpha dS_y=
    \int\limits_{S^k}\omega^\beta dS_y\,.
\end{equation}

  \item[(B)]\label{sphe-pro-inverse-real-Statement} If $\alpha, \beta$ are real, then
\begin{equation}\label{Sphe-Pro-Statement-3}
     \int\limits_{S^k}\omega^\alpha dS_y=
    \int\limits_{S^k}\omega^\beta dS_y\quad \text{implies}\quad \alpha+\beta=k\,\,\,
    \text{or}\,\,\,\alpha=\beta  \,.
\end{equation}

  \item[(C)]\label{sphe-pro-liuville-Statement} For every $\beta\in \mathbb{C}$ there are infinitely many numbers $\alpha\in\mathbb{C}$ such that
\begin{equation}\label{Sphe-Pro-Statement-4}
    \int\limits_{S^k}\omega^\alpha dS_y=
    \int\limits_{S^k}\omega^\beta dS_y\,.
\end{equation}

    \item[(D)]\label{sphe-pro-inverse-complex-Statement} Suppose that for the fixed point $x\notin S^k$
\begin{equation}\label{Sphe-Pro-Statement-5-1}
    \max\{|\Im(\alpha)|, |\Im(\beta)|\}\leq \left.\frac\pi2 \right/ \ln\frac{R+r}{|R-r|}\,.
\end{equation}
Then
\begin{equation}\label{Sphe-Pro-Statement-6}
\int\limits_{S^k} \omega^\alpha dS_y =
\int\limits_{S^k} \omega^\beta dS_y \quad\text{implies}\quad
\alpha+\beta=k\quad\text{or}\quad\alpha=\beta \,.
\end{equation}

\end{description}
\end{theorem}

\begin{remark}
    Statement~(C) shows that the implication introduced in \eqref{Sphe-Pro-Statement-3} fails if we just let $\alpha$ and $\beta$ be complex. However, a certain additional restriction formulated in Statement~(D) for the complex numbers $\alpha, \beta$ allows us to obtain precisely the same implication as in \eqref{Sphe-Pro-Statement-3}.
\end{remark}

\begin{remark}
    Observe that Statement~(B) is a special case of Statement~(D). Indeed, if $\alpha, \beta$ are real, then \eqref{Sphe-Pro-Statement-5-1} holds for every $r\neq R$.
\end{remark}

\begin{corollary}\label{imaginary-part-vanishing-corollary} 
\begin{equation}\label{imaginary-part-vanishing}
    \int\limits_{S^k}\omega^\alpha dS_y=
    \int\limits_{S^k}\omega^{k/2}\cos(b\ln\omega) dS_y\,.
\end{equation}
Moreover, for every natural number $m=0,1,2,3,...$
\begin{equation}\label{logarithm-integral}
\begin{split}
    & \int\limits_{S^k}\omega^{k/2}(\ln\omega)^{2m+1}\cos(b\ln\omega) dS_y=0\,,
    \\& \int\limits_{S^k}\omega^{k/2}(\ln\omega)^{2m}\sin(b\ln\omega) dS_y=0
\end{split}
\end{equation}
and for every $\alpha\in\mathbb{C}$
\begin{equation}\label{Taylor-Decomposition}
    F(\alpha)=\int\limits_{S^k}\omega^{\alpha}dS_y=
    \sum\limits_{m=0}^{\infty}\frac{(\alpha-k/2)^{2m}}{(2m)!}
    \int\limits_{S^k}\omega^{k/2}(\ln\omega)^{2m}dS_y\,,
\end{equation}
which is the Taylor decomposition of $F(\alpha)$ at $\alpha=k/2$.
\end{corollary}

\begin{corollary}\label{sphe-pro-distance-basic-Statement}
If $\alpha+\beta=2k$ and $\alpha,\beta\in\mathbb{C}$, then
\begin{equation}\label{Sphe-Pro-Statement-9.0}
\int\limits_{S^k} \frac{dS_y}{|x-y|^\alpha}
=|R^2-|x|^2|^{(\beta-\alpha)/2}\, \int\limits_{S^k}
\frac{dS_y}{|x-y|^\beta}\,.
\end{equation}
\end{corollary}

\begin{remark} It is clear that formula \eqref{Sphe-Pro-Statement-9.0} can be obtained from Statement~(A) if we rewrite formula \eqref{Sphe-Pro-Statement-2} in terms of distance between $x$ and $y$. In particular, for $\alpha=k-1$ and for every $x\notin S^k$, formula~\eqref{Sphe-Pro-Statement-9.0} yields
\begin{equation}\label{Sphe-Pro-Statement-9}
    \int\limits_{S^k} \frac{dS_y}{|x-y|^{k-1}}
    = \int\limits_{S^k}
    \frac{|R^2-|x|^2|}{|x-y|^{k+1}}\,dS_y\,.
\end{equation}
Note that for $k>1$ the integrand on the left is the Newtonian
potential in $\mathbb{R}^{k+1}$ and the integrand on the right is
exactly the Poisson kernel in $\mathbb{R}^{k+1}$.
\end{remark}

\begin{corollary}\label{One-dimension-sphe-pro} If $k=1$, then for every $p\in\mathbb{C},a,b
\in\mathbb{R}$ and $a>b>0$,
\begin{equation}\label{Sphe-Pro-Statement-10}
    \int\limits_0^{2\pi} \frac{d\theta}{(a-b\,\sin(\theta))^{p}}=
    (a^2-b^2)^{1/2-p}\, \int\limits_0^{2\pi}
    (a-b\,\sin(\theta))^{p-1}d\theta\,.
\end{equation}
\end{corollary}

\begin{corollary}\label{Basic-Th-Integrable-Case-general-Item}
If $k=2$, then
\begin{equation}\label{Basic-Th-Integrable-Case-1}
    \int\limits_{S^2(R)}\omega^{1+ib}dS_y=\frac{2\pi R}{r}
    \cdot\frac{|R^2-r^2|}{b}\cdot\sin\left(b\ln\frac{R+r}{|R-r|}\right)\,.
\end{equation}
In particular, if $b=0$,
\begin{equation}\label{Basic-Th-Integrable-Case-2}
    \int\limits_{S^2(R)}\omega\, dS_y=\frac{2\pi R}{r}\cdot|R^2-r^2|\cdot\ln\frac{R+r}{|R-r|}\,.
\end{equation}
\end{corollary}

\section{Proof of the basic results}

\begin{proof}[Proof of Theorem \ref{Basic-theorem}.] We are going to prove successively  all of the Statements (A), (B), (C) and (D) from the Theorem.

\begin{proof}[Proof of Statement (A), p.~\pageref{Sphe-pro-direct-Statement}.]

\textbf{Case 1 ($|x|<R$).} For
notation refer to Figure~\ref{Change}. Let $y$ be an arbitrary
point of $\mathbb{S}^k$. Here, $y^\star$ is defined as an
intersection of the sphere $\mathbb{S}^k$ and line $xy$. Angle
$\theta$ is defined as $\pi - \angle xOy$ and $\theta^\star$ is
defined as $\pi-\angle xOy^\star$. The circle $\mathbb{S}^1$
appears on the cross-sectional plane defined by the
points $O$, $x$ and $y\in \mathbb{S}^k$. Let $n(\theta)=|xy|$ and
$m(\theta)=|y^\star x|$ and let the segments $a$ be
perpendicular to $xO$.

\begin{figure}[h]
    \center\epsfig{figure=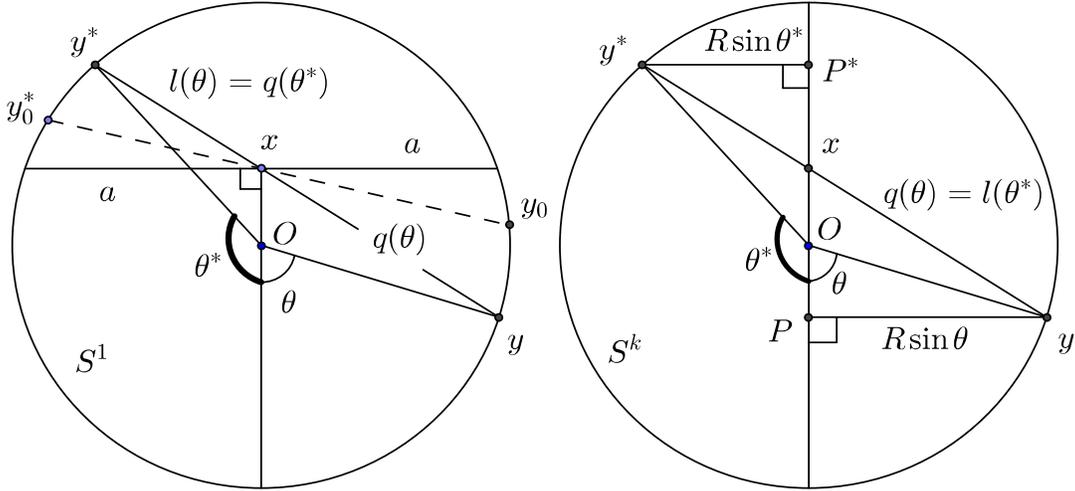, height=6.7cm, width=14.5cm}
    \caption{Change of variables}
    \label{Change}
\end{figure}



One can easily check that
\begin{equation}\label{6}
    \omega=\frac{R^2-|x|^2}{|x-y|^2}=\frac{|a|^2}{q^2}=\frac{q\cdot
    l}{q^2}=\frac{l}{q}\,\, ,
\end{equation}
since $|a|^2=q\cdot l$ and, by Pythagorean theorem,
$|a|^2=R^2-|x|^2$.

In the next step we are going to use the change of variables  from $\theta$ to $\theta^*$, introduced also by Hermann Schwarz in~\cite{Schwarz} (pp. 359-361).
Notice, first, that $\triangle y_0xy$ is similar to $\triangle
xy_0y$ and therefore,
\begin{equation}\label{2-z}
    \frac{d\theta^\star}{d\theta}=-\lim_{y_0\rightarrow y}
    \frac{|y^\star_0 y^\star|}{|y_0 y|} =
    -\frac{l(\theta)}{q(\theta)}\,.
\end{equation}
It is not hard to see that
\begin{equation}\label{7-z}
    l(\theta)=|y^\star x|=q(\theta^\star) \quad
    \text{and} \quad
    q(\theta)=|xy|=l(\theta^\star)\,.
\end{equation}
Therefore, \eqref{2-z} combined with \eqref{7-z} yields
\begin{equation}\label{3-z}
    d\theta=-\frac{q(\theta)}{l(\theta)}d\theta^\star=
    -\frac{l(\theta^\star)}{q(\theta^\star)} d\theta^\star\,.
\end{equation}
Let $P$ and $P^\star$ be orthogonal projections of $y$ and
$y^\star$ to line $xO$ respectively. Clearly, $\triangle
xy^\star P^\star$ is similar to $\triangle xyP$ and thus
\begin{equation}\label{4-z}
    l(\theta)\,
    R \sin\theta=q(\theta)\,
    R\sin\theta^\star=
    l(\theta^\star)
    \, R\sin\theta^\star\,.
\end{equation}
Using \eqref{6}, \eqref{7-z}, \eqref{3-z} and \eqref{4-z}, we have

\begin{equation}\label{5-z}
\begin{split}
  & \oint\limits_{\mathbb{S}^k} \omega^\alpha d\,\mathbb{S}_y=
    \int\limits_0^\pi
    \left(\frac{l(\theta)}{q(\theta)}\right)^\alpha
    (R\sin\theta))^{k-1} |\sigma_{k-1}| Rd\theta
    \\& =
    \int\limits_0^\pi
    \left(\frac{l(\theta^\star)}{q(\theta^\star)}\right)^{k-\alpha}(R\sin\theta^\star)^{k-1}
    |\sigma_{k-1}|Rd\theta^\star
    = \oint\limits_{\mathbb{S}^k}
    \omega^{k-\alpha} d\,\mathbb{S}_y\,,
\end{split}
\end{equation}
where $|\sigma_{k-1}|$ is the area of $(k-1)$-dimensional unit
sphere. The formula
$d\,\mathbb{S}_y=(R\sin\theta)^{k-1}|\sigma_{k-1}|R\,d\theta$
holds, since our integrand $\omega(x,y)$ depends only on angle
$\theta$ when $x$ is fixed. 
$\blacksquare$ \medskip

\textbf{Case 2 ($|x|>R$).} Let
$\overline{x}\in Ox$ 
and $|\overline{x}|\cdot |x|=R\,^2$. Then, one can easily check
that $\omega(x,y)=\omega(\overline{x},y)$ for every
$y\in\mathbb{S}^k$. Clearly, the point $\overline{x}$ is inside
the sphere $\mathbb{S}^k$ and then we can apply the result obtained in formula \eqref{5-z}. This yields
\begin{equation}\label{1-k}
    \oint\limits_{\mathbb{S}^k} \omega^{\alpha}(x,y) d\,\mathbb{S}_y=
    \oint\limits_{\mathbb{S}^k}
    \omega^{\alpha}(\overline{x},y) d\,\mathbb{S}_x=
    \oint\limits_{\mathbb{S}^k}
    \omega^{\beta}(\overline{x},y) d\,\mathbb{S}_x = \oint\limits_{\mathbb{S}^k}
    \omega^{\beta}(x,y) d\,\mathbb{S}_y\,.
\end{equation}
This completes the proof of Statement~(A) from p.~\pageref{Sphe-pro-direct-Statement}.
\end{proof}

\begin{proof}[Proof of Statement~(B), p.\pageref{sphe-pro-inverse-real-Statement}.]

Fix a real number $\alpha$. We have to show that the equation
\begin{equation}\label{2-r}
    F(\alpha)=\oint\limits_{\mathbb{S}^k} \omega^{\alpha}d\,\mathbb{S}_y=\oint\limits_{\mathbb{S}^k}
    \omega^{\lambda}d\,\mathbb{S}_y=F(\lambda)
\end{equation}
has only two solutions with respect to the variable $\lambda$;
namely, $\lambda=\alpha$ and $\lambda=k-\alpha$. First note that,
since $\lambda+(k-\lambda)=k$, we have $F(\lambda)=F(k-\lambda)$
for every $\lambda\in \mathbb{R}$. Therefore, the function
$F(\lambda)$ is symmetric with respect to the point $\lambda=k/2$.
The Leibnitz rule of differentiation
under the integral sign yields

\begin{equation}\label{3-r}
    \frac{dF(\lambda)}{d\lambda}=\oint\limits_{\mathbb{S}^k}
    \omega^{\lambda}\, (\ln\omega)\,d\,\mathbb{S}_y\,.
\end{equation}
We can apply the same rule one more time,
\begin{equation}\label{4-r}
    \frac{d^2F(\lambda)}{(d\lambda)^2}=\oint\limits_{\mathbb{S}^k}
    \omega^{\lambda}\, (\ln\omega)^2\,d\,\mathbb{S}_y >0\quad
    \forall\lambda\in \mathbb{R}\,.
\end{equation}
Therefore, the graph of our function $\xi=F(\lambda)$, plotted in
the $(\lambda,\xi)$-plane, is convex and symmetric with respect to
the line $\lambda=k/2$. This means that $\lambda=k/2$ is the
minimum point for $F(\lambda)$. Therefore, for every $\nu >
F(k/2)$, the equation $F(\lambda)=\nu$ has only two solutions
$\lambda=\alpha$ and $\lambda=k-\alpha$. This completes the proof of Statement~(B) from p.\pageref{sphe-pro-inverse-real-Statement}.
\end{proof}

\begin{proof}[Proof of Statement~(C), p.~\pageref{sphe-pro-liuville-Statement}]
Let $F(\alpha)=\int\limits_{S^k} \omega^\alpha dS_y$, where $\alpha=u+iv$. Then, the Statement~(C) can be obtained as a consequence of the following Picard's Great Theorem.

\begin{theorem}[Picard's Great Theorem]\label{Picard-Great-theorem}
 Let $c\in\mathbf{C}$ be an isolated essential singularity of $f$. Then, in every neighborhood of $c$, $f$ assumes every complex number as a value with at most one exception infinitely many times, see \cite{Reinhold}, (p.240).
\end{theorem}

Indeed, the direct computation shows that $F$ is an entire function, which means that $\infty$ is an isolated singularity of $F$.
Note that the Taylor decomposition presented in \eqref{Taylor-Decomposition}, p.~\pageref{Taylor-Decomposition} shows that $\infty$ is the essential singularity for $F(\alpha)$, since all even Taylor coefficients are positive.

It is clear also, we can not meet an exception mentioned in Picard's theorem above, since we are looking for complex numbers $\lambda$, satisfying $F(\lambda)=F(\alpha)$ for some $\alpha$. This means that the value $F(\alpha)$ is already assumed by $F$, and then, by Picard's Great Theorem, $F$ must attain this value infinitely many times in every neighborhood of an isolated essentially singular point, i.e., at every neighborhood of $\infty$, in our case. This completes the proof of Statement~(C).
\end{proof}

\begin{proof}[\textbf{Proof of Statement (D), p.~\pageref{sphe-pro-inverse-complex-Statement}.}]

Denote
\begin{equation}\label{E-basic-proof-1}
    \Upsilon=\Upsilon(x)=\min\left\{\frac{R}{|x|},\frac{|x|}{R}\right\}.
\end{equation}
Recall that $R$ and $x\notin S^k(R)$ are fixed. Then,
\begin{equation}\label{E-basic-proof-2}
    \omega(x,y)=\omega(\theta)=\frac{1-\Upsilon^2}{1+\Upsilon^2+2\Upsilon\cos\theta}\,,
\end{equation}
where $\theta=\pi-\angle xOy$, $\theta\in[0,\pi]$ and it is clear that $\Upsilon<1$ for every possible $x\notin S^k(R)$. Observe also that $\omega(\theta)$ is increasing while $\theta\in[0,\pi]$. Therefore,
\begin{equation}\label{E-basic-proof-3}
    \frac{1-\Upsilon}{1+\Upsilon}\leq\omega(\theta)\leq\frac{1+\Upsilon}{1-\Upsilon}.
\end{equation}
Let $p\geq 0$ be an arbitrary non-negative number such that
\begin{equation}\label{introduction-of-p}
    \max\{|\Im(\alpha)|, |\Im(\beta)|\}\leq p \leq \left.\frac\pi2\right/\ln\frac{R+r}{|R-r|}=\left.\frac\pi2\right/\ln\frac{1+\Upsilon}{1-\Upsilon}\,.
\end{equation}
The direct computations show that the second inequality in~\eqref{introduction-of-p}, in itself, implies
\begin{equation}\label{E-basic-proof-4}
    e^{-\pi/2p}\leq\frac{1-\Upsilon}{1+\Upsilon}\quad\text{and}\quad
    e^{\pi/2p}\geq\frac{1+\Upsilon}{1-\Upsilon}
\end{equation}
Therefore, combining \eqref{E-basic-proof-3} and \eqref{E-basic-proof-4}, we have
\begin{equation}\label{E-basic-proof-5}
    e^{-\pi/2p}\leq\frac{1-\Upsilon}{1+\Upsilon}\leq\omega(\theta)\leq
    \frac{1+\Upsilon}{1-\Upsilon}\leq e^{\pi/2p}\,,
\end{equation}
or, equivalently,
\begin{equation}\label{E-basic-proof-6}
    -\frac{\pi}{2}\leq p\ln\omega\leq\frac{\pi}{2}\quad\text{or}\quad
    p\ln\omega\in\left[-\frac{\pi}{2},\frac{\pi}{2}\right]\,.
\end{equation}
If $\alpha=\xi+i\zeta$ and $|\zeta|\leq p$, then, clearly, \eqref{E-basic-proof-6} implies
\begin{equation}\label{E-basic-proof-7}
    -\frac{\pi}{2}\leq \zeta\ln\omega\leq\frac{\pi}{2}\quad\text{or}\quad
    \zeta\ln\omega\in\left[-\frac{\pi}{2},\frac{\pi}{2}\right]\,.
\end{equation}
Hence,
\begin{equation}\label{cosine-positivenes}
\cos(\zeta\ln\omega(\theta))\geq 0 \quad\text{for every}\quad \zeta\in[-p,p]
\end{equation}
and for every $\theta\in[0,\pi]$. Fix a number $p>0$ satisfying \eqref{cosine-positivenes} or, equivalently~\eqref{Sphe-Pro-Statement-5-1}, p.~\pageref{Sphe-Pro-Statement-5-1} and introduce the set
\begin{equation}\label{Strip-definition}
    \Xi=\{\xi+i\zeta\mid-p\leq\zeta\leq p\}\,.
\end{equation}
Consider the function
\begin{equation}\label{Patch-6}
    F(\alpha)=\int\limits_{S^k} \omega^\xi\cos(\zeta\ln\omega) dS_y +
    i\int\limits_{S^k} \omega^\xi\sin(\zeta\ln\omega) dS_y \,.
\end{equation}
Note now that Statement~(E) will be proven if we show that for any $\beta=a+ib\in\Xi$, the equation $F(\alpha)=F(\beta)$ has only two solutions in $\Xi$, i.e., $\alpha=\beta$ and $\alpha=k-\beta$. To accomplish this goal let us split our argument into the following four steps.
\begin{description}
  \item[Step 1.] Fix any $\beta=a+ib\in\Xi$ and set the equation $F(\alpha)=F(\beta)$ as the following system of two equations
      \begin{equation}\label{system_for_real_Imaginary}
        \Re F(\xi,\zeta)=\Re F(\beta)\quad\text{and}\quad \Im F(\xi,\zeta)=\Im F(\beta).
      \end{equation}
           Then, we establish a symmetry property for $\Re F(\xi,\zeta)$ and $\Im F(\xi,\zeta)$. We shall see that $\Re F(\xi,\zeta)$ is symmetric with respect to the two lines $\xi=k/2$ and $\zeta=0$, while $\Im F(\xi,\zeta)$ is skew-symmetric with respect to the same lines. This step is described on p.~\pageref{Step_curve_1}.
  \item[Step 2.] Here we describe all solutions inside $\Xi$ for the equation $\Re F(\xi,\zeta)=\Re F(\beta)$. Because of the symmetry mentioned in Step 1 it is enough to describe the set of solutions in the first quadrant $[k/2,\infty)\times[0,p]$. The description will be presented in Proposition~\ref{Level-curve-description}, page~\pageref{Level-curve-description}.


  \item[Step 3.] Here we shall study the behavior of $\Im F(\xi,\zeta)$ along any of the $\mathcal{C}^1$ curves from Step 2. Proposition~\ref{uniquiness-total-in-first-quadrant}, p.~\pageref{uniquiness-total-in-first-quadrant} shows that $\Im F(\xi,\zeta)$ is a strictly increasing function along any of the level curves inside the first quadrant $[k/2,\infty)\times[0,p]$, which implies that the system \eqref{system_for_real_Imaginary} has the unique solution in the first quadrant for any $\beta\in[k/2,\infty)\times[0,p]$, i.e., the solution is $\beta$ itself. This step is carried out on p.~\pageref{Step_curve_3}.
  \item[Step 4.] Here we study the behavior of $\Im F(\xi,\zeta)$ in the whole strip $\Xi$. We shall see that for every $\beta\in[k/2,\infty)\times[0,p]$ there exist only two solutions of the system \eqref{system_for_real_Imaginary}, i.e., $\alpha_1=\xi+i\zeta=\beta$ and $\alpha_2=k-\beta$. The same argument will show that for every $\beta\in\Xi$ there are at most two possible solutions of system~\eqref{system_for_real_Imaginary}, i.e., $\alpha_1=\beta$ and $\alpha_2=k-\beta$. This step is carried out on p.~\pageref{Step_curve_4}.
\end{description}
Now let us follow this plan described in the Steps 1-4 above.

\underline{\textbf{Step 1.}}\label{Step_curve_1} It is convenient to introduce the notation for the real and the imaginary parts of $F(\alpha)=F(\xi+i\zeta)=F(\xi,\zeta)$. Since $R=|y|$ and $x\notin S^k(R)$ are fixed, we can denote
\begin{equation}\label{Definition-W}
    W(\alpha)=W(\xi+i\zeta)=W(\xi,\zeta)=
    \int\limits_{S^k} \omega^\xi\cos(\zeta\ln\omega) dS_y
\end{equation}
and similarly,
\begin{equation}\label{Definition-I}
    I(\alpha)=I(\xi+i\zeta)=I(\xi,\zeta)=
    \int\limits_{S^k} \omega^\xi\sin(\zeta\ln\omega) dS_y \,.
\end{equation}
Fix some $\beta=a+ib\in\Xi$ and rewrite system \eqref{system_for_real_Imaginary} in the new notation.
\begin{equation}\label{system_for_real_Ima_new_notation}
    W(\xi,\zeta)=W(a,b)\quad\text{and}\quad I(\xi,\zeta)=I(a,b)\,,
\end{equation}
which is going to be solved. As we saw in Statement (A), $F(\alpha)=F(k-\alpha)$ for every $\alpha=\xi+i\zeta\in \mathbb{C}$ and then, the direct computation yields
\begin{equation}\label{Symmetry-Reations}
\begin{split}
    & W(\xi+i\zeta)=W(k-\xi,-\zeta)=W(\xi,-\zeta)=W(k-\xi,\zeta) \quad\text{and}
    \\& I(\xi+i\zeta)=I(k-\xi,-\zeta)=-I(\xi, -\zeta)=-I(k-\xi,\zeta)\,,
\end{split}
\end{equation}
which pictured on the diagram below, see Figure~\ref{Re-Symmetry}.

\begin{figure}[!h]
    \centering
    \epsfig{figure=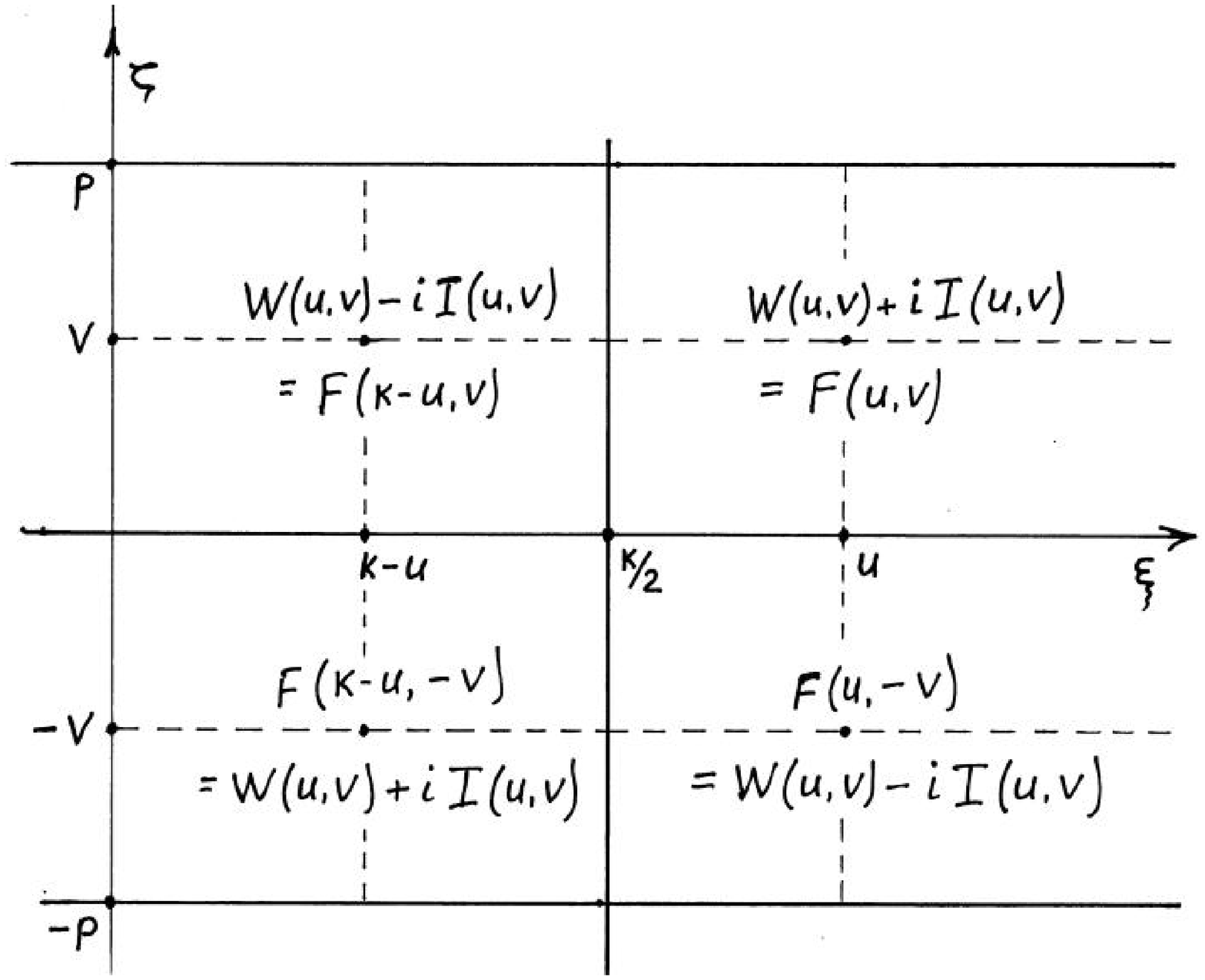,height=6cm,width=9cm}\\
    \caption{The symmetry of $\Re$.}\label{Re-Symmetry}
\end{figure}

The symmetry relations presented in \eqref{Symmetry-Reations} imply that $W(u+iv)$ is symmetric with respect to the two perpendicular lines $\xi=k/2$ and $\zeta=0$ in the plane $(\xi,\zeta)$. Similarly, $I(\xi,\zeta)$ is symmetric with respect to point $(k/2,0)$ and skew-symmetric with respect to the lines $\zeta=0$ and $\xi=k/2$.

\underline{\textbf{Step 2.}}\label{Step_curve_2} Because of the symmetry, without loss of generality, we may assume that $\beta=a+ib\in\{[k/2,\infty)\times[0,p]\}$, which is the upper-right quadrant. Now we are ready to solve the equation $F(\xi+i\zeta)=F(\beta)$. First we need to describe all solutions for the real parts of this equation, i.e., $\Re F(\xi,\zeta)=\Re F(\beta)$.
Let us define two quadrants $Q$, $\overline{Q}$ and the level set $S(a,b)$ by
\begin{equation}\label{Definition-Q-S(a,b)}
\begin{split}
    & Q=(k/2,\infty)\times(0,p)\,,\quad \overline{Q}=[k/2, \infty)\times[0,p],
    \\& S(a,b)=\{(\xi,\zeta)\in\overline{Q}\mid W(\xi,\zeta)=W(a,b)\}\,.
\end{split}
\end{equation}
It is clear that $\overline{Q}$ is the closure of $Q$.
The following proposition describes all possible types of such level sets.

\begin{proposition}[Level curves description]\label{Level-curve-description} \
    \begin{enumerate}
      \item For every fixed $\beta=a+ib\in\overline{Q}$ the set of solutions $S(a,b)$ for the equation
        \begin{equation}
            W(\xi,\zeta)=\int\limits_{S^k}\omega^\xi\cos(\zeta\ln\omega)dS_y=W(a,b)
        \end{equation}
        is a continuous curve $\gamma_{(a,b)}(t)=\gamma(t)=(t,v(t))\subseteq{\overline{Q}}$.
      \item The function $v(t)$ is $\mathcal{C}^{\infty}$ and $v'(t)>0$ for every point $t$ such that $(t,v(t))\subseteq{Q}$.
      \item No two of these curves have common points in $\overline{Q}$.
    \end{enumerate}
\end{proposition}

\begin{remark}
    According to Theorem 6.1, (\cite{Artamoshin_Diss},~p.~128), if a level curve starts at the corner $C=(k/2,0)$, then it bisects the corner and if a level curve starts at the lower or at the left edge of $\overline{Q}$, then it must be perpendicular to the edge. No two of these level curves have a common point in $\overline{Q}$.
\end{remark}

\begin{proof}[\textbf{Proof of the Proposition \ref{Level-curve-description}}]

First, we need to analyze the behavior of partial derivatives of $W(\xi,\zeta)$ in $\overline{Q}$.
\begin{claim}\label{Partial-dW-d-xi-claim}
    \begin{equation}\label{Partial-dW-d-xi}
    \begin{split}
        & \frac{\partial W(\xi,\zeta)}{\partial\xi}>0\quad\text{for every}\quad(\xi,\zeta)\in
        \left(\frac{k}{2},\infty\right)\times[0,p]\,;
        \\&  \frac{\partial W(\xi,\zeta)}{\partial\xi}=0\quad\text{if}\quad\xi=
        \frac{k}{2}\quad\text{and}\quad\zeta\in[0,p]\,.
    \end{split}
    \end{equation}
\end{claim}

\begin{proof}[Proof of the claim \ref{Partial-dW-d-xi-claim}]
    Note that $W(\xi,\zeta)$ is $\mathcal{C}^{\infty}$ function of both variables and, according to \eqref{Symmetry-Reations}, $W(\xi,\zeta)$ is symmetric with respect to $\xi=k/2$ as a function of $\xi$ for every $\zeta$. This implies that
    \begin{equation}\label{horizontal-derivative-zero}
        \frac{\partial W(\xi,\zeta)}{\partial\xi}=0\quad\text{for}\quad\xi=
        \frac{k}{2}\quad\text{and for every}\quad\zeta\,.
    \end{equation}
    Note also that according to \eqref{cosine-positivenes}, p.~\pageref{cosine-positivenes},
    \begin{equation}\label{horizontal-second-derivative}
        \frac{\partial^2 W(\xi,\zeta)}{\partial\xi^2}=
        \int\limits_{S^k}\omega^{\xi}(\ln\omega)^2\cos(\zeta\ln\omega)dS_y>0
    \end{equation}
    for every $(\xi,\zeta)\in\Xi$ defined in \eqref{Strip-definition}. It is clear that \eqref{horizontal-derivative-zero} and \eqref{horizontal-second-derivative} imply \eqref{Partial-dW-d-xi}. This completes the proof of Claim~\ref{Partial-dW-d-xi-claim}.
    \end{proof}

    \begin{claim}\label{partial-vertical-claim}
      \begin{equation}\label{partial-vertical}
      \begin{split}
        & \frac{\partial W(\xi,\zeta)}{\partial\zeta}<0\quad\text{for every}\quad(\xi,\zeta)\in
        \left[\frac{k}{2},\infty\right)\times(0,p]\,;
        \\&  \frac{\partial W(\xi,\zeta)}{\partial\zeta}=0\quad\text{if}\quad\zeta=0
        \quad\text{and}\quad\xi\in[k/2, \infty]\,.
      \end{split}
      \end{equation}
    \end{claim}
    \begin{proof}[Proof of the Claim \ref{partial-vertical-claim}]
        Here we are going to use similar symmetry and second derivative arguments used in the proof of previous Claim. Recall that $W(\xi,\zeta)$ is $\mathcal{C}^{\infty}$ function of both variables and, according to \eqref{Symmetry-Reations}, $W(\xi,\zeta)$ is symmetric with respect to $\zeta=0$ as a function of $\zeta$ for every $\xi$. This implies that
    \begin{equation}\label{vertical-derivative-zero}
        \frac{\partial W(\xi,\zeta)}{\partial\zeta}=0\quad\text{for}\quad\zeta=0
        \quad\text{and for every}\quad\xi\,.
    \end{equation}
    Note also that according to \eqref{cosine-positivenes},
    \begin{equation}\label{vertical-second-derivative}
        \frac{\partial^2 W(\xi,\zeta)}{\partial\zeta^2}=
        -\int\limits_{S^k}\omega^{\xi}(\ln\omega)^2\cos(\zeta\ln\omega)dS_y<0
    \end{equation}
    for every $(\xi,\zeta)\in\Xi$ defined in \eqref{Strip-definition}, p.~\pageref{Strip-definition}. It is clear that \eqref{vertical-derivative-zero} and \eqref{vertical-second-derivative} imply \eqref{partial-vertical}. This completes the proof of Claim~\ref{partial-vertical-claim}.
    \end{proof}

    For the next step we need the following notation for the boundary parts of $\overline{Q}$. With every $\beta=a+ib\in\overline{Q}$ we associate two curves $\Gamma_0(a,b)$ and $\Gamma_1(a,b)$ as shown on the Figure~\ref{Special-boundary-curves} below and defined as follows.
\begin{equation}
\begin{split}
     \Gamma_0(a,b)= & \{(\xi,\zeta)\mid\xi=k/2\,\,\text{and}\,\,\zeta\in[0,b]\}\cup
    \\& \{(\xi,\zeta)\mid\xi\in[k/2,a]\,\,\text{and}\,\,\zeta=0\};
\end{split}
\end{equation}
\begin{equation}
    \Gamma_1(a,b)=\{(\xi,\zeta)\mid \xi\in[a,\infty)\,\,\text{and}\,\,\zeta=p\}\,.
\end{equation}

\begin{figure}[!h]
    \centering
    \epsfig{figure=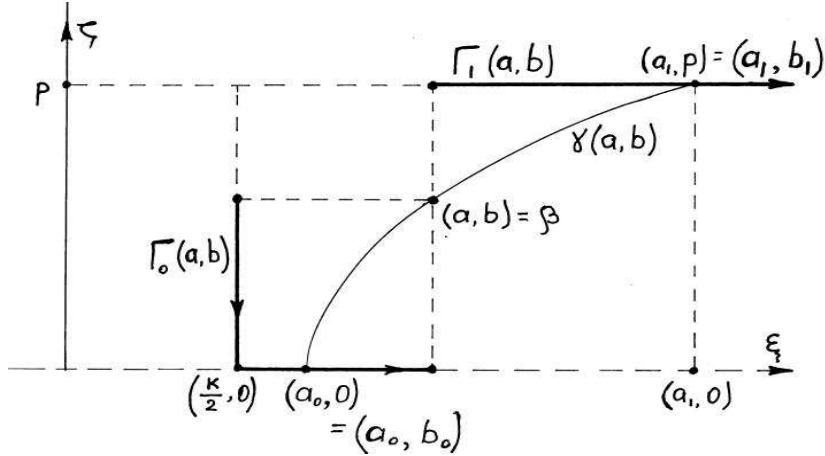,height=6cm,width=11cm}
  \caption{Special boundary curves.}\label{Special-boundary-curves}
\end{figure}

Note then, according to \eqref{Partial-dW-d-xi} and \eqref{partial-vertical}, the value of $W(\xi,\zeta)$ is a continuous and strictly increasing function as the point $(\xi, \zeta)$ runs along $\Gamma_0(a,b)$ from $(k/2,b)$ to $(a,0)$ and
\begin{equation}
    W(k/2,b)\leq W(\xi,\zeta)\leq W(a,0)\quad\text{for every}\,\,(\xi,\zeta)\in\Gamma_0(a,b)\,.
\end{equation}
On the other hand, $W(k/2,b)\leq W(a,b)\leq W(a,0)$.
Therefore, there exists the unique point $(a_0,b_0)\in\Gamma_0$ such that $W(a_0,b_0)=W(a,b)$.

Note also that by \eqref{partial-vertical}, $W(\xi,\zeta)$ is continuous and strictly decreasing along the vertical segment connecting $(a,b)$ and $(a,p)$. Thus, $W(a,p)\leq W(a,b)$. Then, when the point $(\xi,\zeta)$ runs along $\Gamma_1$ from $(a,p)$ to $\infty$, the value of $W(\xi,\zeta)$ is continuously and strictly increasing to $\infty$, according to~\eqref{Partial-dW-d-xi} and \eqref{horizontal-second-derivative}. Therefore, there exists the unique point $(a_1,p)\in\Gamma_1$ such that $W(a_1,p)=W(a,b)$. Clearly, $a_0\leq a\leq a_1$.

Now we are ready to describe the level set $S(a,b)$ for any fixed $(a,b)\in \overline{Q}$.

\begin{claim}\label{Level-set-boundary}
    \begin{equation}\label{rectangle-for-level-set}
        S(a,b)\subseteq[a_0,a_1]\times[0,p]\,.
    \end{equation}
\end{claim}

\begin{proof}[Proof of the Claim \ref{Level-set-boundary}]

Using \eqref{Partial-dW-d-xi} and \eqref{partial-vertical} we can observe that the value of $W(a_0,0)$ is the strict maximum of $W(\xi,\zeta)$ for $(\xi,\zeta)\in [k/2,a_0]\times[0,p]$ as well as the value of $W(a,p)$ is the strict minimum for $(\xi,\zeta)\in [a_1,\infty]\times[0,p]$. This completes the proof of Claim~\ref{Level-set-boundary}.
\end{proof}

\begin{claim}\label{function-presentation-claim} $S(a,b)$ can be interpreted as the graph of a function $v(t)$, i.e.,
\begin{equation}\label{function-presentation}
    S(a,b)=\{(t,v(t))\mid t\in[a_0,a_1]\}\,.
\end{equation}
\end{claim}

\begin{proof} Fix any $t\in[a_0,a_1]$. Then, for the vertical segment connecting two points $(t,0)$ and $(t,p)$ we can observe that
\begin{equation}
    W(t,0)\geq W(a,b)\geq W(t,p)
\end{equation}
since, according to \eqref{Partial-dW-d-xi}, the value of $W(\xi,\zeta)$ is a strictly increasing function along every horizontal line in $\overline{Q}$. Note also that the value of $W(t,\zeta)$, according to \eqref{partial-vertical}, is strictly decreasing for $\zeta\in[0,p]$. Hence, there exists the unique value of $\zeta=v(t)\in[0,p]$, such that $W(t,v(t))=W(a,b)$, which completes the proof of Claim \ref{function-presentation-claim}.
\end{proof}

To complete the proof of Proposition~\ref{Level-curve-description}, from p.~\pageref{Level-curve-description}, we need to study the functional properties of the function $v(t)$. First, let us observe that $v(t)>0$ for every $t\in(a_0,a_1)$, since the value of $W(\xi,\zeta)$, according to \eqref{Partial-dW-d-xi}, is a strictly increasing function along the lower edge of $\overline{Q}$. Therefore, by the Implicit Function Theorem, for every $t\in(a_0,a_1)$ there exists some neighborhood $N_n(t)=(t-\epsilon_n, t+\epsilon_n)$, where $v(t)$ is continuously differentiable $n$ times for any natural $n$ and then, according to \eqref{Partial-dW-d-xi} and \eqref{partial-vertical}, we have
      \begin{equation}\label{derivative-v-by-Implicit-FT}
        \frac{dv(t)}{dt}=\frac{W_\xi(t,v(t))}{-W_\zeta(t,v(t))}>0\quad\text{for every}\,\,t\in(a_0,a_1)\,.
      \end{equation}
      A possible curve $\gamma(t)=(t,v(t))$ starting with the horizontal edge of $\overline{Q}$ is sketched on Figure \ref{Special-boundary-curves}, page~\pageref{Special-boundary-curves}.
Note also that $v(t)$ is continuous in $[a_0,a_1]$, since $v(t)$ is monotone in $(a_0,a_1)$ and $W(\xi,\zeta)$ is continuous in $\overline{Q}$. Now the proof of the Proposition \ref{Level-curve-description} from p.~\pageref{Level-curve-description} is complete.
\end{proof}

\underline{\textbf{Step 3.}}\label{Step_curve_3} The goal of this step is to obtain the following proposition.

\begin{proposition}\label{uniquiness-total-in-first-quadrant}
   If $(a,b)\in\overline{Q}$, then the equation $F(\xi,\zeta)=F(a,b)$ has the unique solution $(\xi,\zeta)\in\overline{Q}$, i.e., $(\xi,\zeta)=(a,b)$.
\end{proposition}

       \begin{proof}[Proof of the Corollary \ref{uniquiness-total-in-first-quadrant}.]
        It is clear that a solution $(\xi,\zeta)$ may appear only on the level curve $S(a,b)=\gamma(t)$ described above. For a point $(\xi,\zeta)\in S(a,b)$ to be a solution we need $I(\xi,\zeta)=I(t,v(t))=I(a,b)$. So, let us study the behavior of $I(\xi,\zeta)$ along $\gamma(t)$. By the direct computation we can observe that
        \begin{equation}\label{Cauchy-Riemann}
            \frac{\partial W}{\partial \xi}=\frac{\partial I}{\partial\zeta}\quad\text{and}\quad
            \frac{\partial W}{\partial \zeta}=-\frac{\partial I}{\partial\xi}\,,
        \end{equation}
        which together with \eqref{derivative-v-by-Implicit-FT} leads to
        \begin{equation}
            \frac{d I(t,v(t))}{dt}=
            \left( 1+[v'(t)]^2\right)\left[-\frac{\partial W(t,v(t))}{\partial v(t)}\right]>0
        \end{equation}
        for every $t\in(a_0,a_1)$, where the last inequality holds because of \eqref{partial-vertical} from Claim~\ref{partial-vertical-claim}, p.~\pageref{partial-vertical-claim}. Therefore, the value of $I(t,v(t))$ is a strictly increasing function on $[a_0, a_1]$ and then, the function $I(t,v(t))$ assumes each of its values for $t\in[a_0,a_1]$ only once. Hence, the equation $I(\xi,\zeta)=I(a,b)$ must have only one solution in $S(a,b)$. This completes the proof of the Proposition \ref{uniquiness-total-in-first-quadrant}.
      \end{proof}

\underline{\textbf{Step 4.}}\label{Step_curve_4} The next and the last stage in the proof of Statement (D) is to analyze the behavior of $I(\xi,\zeta)=\Im F(\alpha)$ in $\Xi=\{(\xi,\zeta)\mid\zeta\in[-p,p]\}.$

\begin{proposition}\label{Signum-Behavior-of-I-propo}
    \begin{equation}\label{Signum-Behavior-of-I-1}
        I(\xi,\zeta)>0\quad\text{for}\quad(\xi,\zeta)\in \{(0,\infty)\times(0,p]\}\cup\{(-\infty,0)\times(0,-p]\};
    \end{equation}

    \begin{equation}\label{Signum-Behavior-of-I-2}
        I(\xi,\zeta)<0\quad\text{for}\quad(\xi,\zeta)\in \{(-\infty,0)\times(0,p]\}\cup\{(0,\infty)\times(0,-p]\}.
    \end{equation}
    The signum of $I$ behavior is sketched on Figure~\ref{Signum-Behavior-of-I} below.

  \begin{figure}[!h]
    \centering
    \epsfig{figure=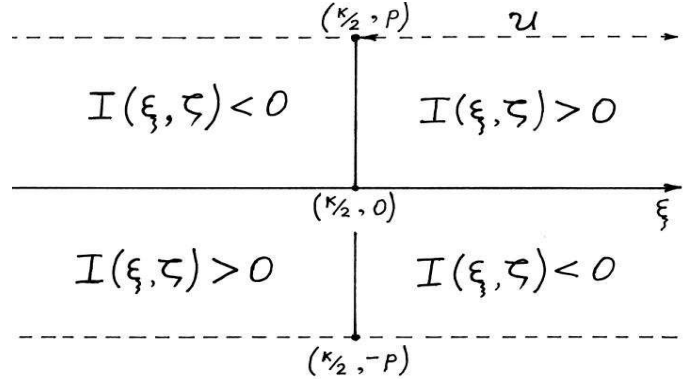,height=5cm,width=9cm}\\
    \caption{Signum behavior of $I(\xi,\zeta)$.}\label{Signum-Behavior-of-I}
\end{figure}
\end{proposition}

\begin{proof}[Proof of Proposition \ref{Signum-Behavior-of-I-propo}.]
    Using the skew-symmetry of $I(\xi,\zeta)$ introduced in \eqref{Symmetry-Reations}, p.~\pageref{Symmetry-Reations}, we can observe that
\begin{equation}\label{cross-imaginary-vanishing}
    I(\xi,\zeta)=\int\limits_{S^k}\omega^\xi\sin(\zeta\ln\omega)dS_y=0\quad\text{for}\quad
    \xi=0\,\,\text{or}\,\,\zeta=k/2\,.
\end{equation}
On the other hand, we can observe that the condition $\zeta\ln\omega\in[-\pi/2,\pi/2]$ yields
\begin{equation}\label{signum-relation-for-sinus}
    \text{sign}(\ln\omega\cdot\sin(\zeta\ln\omega))=\text{sign}\zeta
\end{equation}
for every value of $\omega$, except $\omega=1$ and then, for every point $(\xi,\zeta)\in\Xi$, one of following relations depending on the signum of $\zeta$ holds
\begin{equation}\label{Cross-difinition}
\begin{split}
    & \left.\frac{dI(\xi,\zeta)}{d\xi}\right|_{\zeta=0}=
    \left[\int\limits_{S^k}\omega^\xi(\ln\omega)\sin(\zeta\ln\omega)dS_y \right]_{\zeta=0}=0;
    \\& \left.\frac{dI(\xi,\zeta)}{d\xi}\right|_{\zeta>0}>0 \quad\text{or}\quad\left.\frac{dI(\xi,\zeta)}{d\xi}\right|_{\zeta<0}<0,
\end{split}
\end{equation}
which to together with \eqref{cross-imaginary-vanishing} implies \eqref{Signum-Behavior-of-I-1} and \eqref{Signum-Behavior-of-I-2}. This completes the proof of Proposition~\ref{Signum-Behavior-of-I-propo}.
\end{proof}


Note also that $I(\xi,\zeta)$, as well as $W(\xi,\zeta)$ are symmetric with respect to $(k/2,0)$. Therefore, the uniqueness obtained in Proposition~\ref{uniquiness-total-in-first-quadrant}, p.~\pageref{uniquiness-total-in-first-quadrant} together with the symmetry of $W$ and the skew-symmetry of $I$ with respect to the lines $\xi=k/2$ and $\zeta=0$ implies that the equation $F(\xi,\zeta)=F(a,b)$ has only two solutions $(\xi,\zeta)=(a,b)$ and $(\xi,\zeta)=(k-a,-b)$ for every $(a,b)\in\overline{Q}$.

It is clear that a similar symmetry argument shows that for every $(a,b)\in\Xi$, the equation $F(\xi,\zeta)=F(a,b)$ also must have only two solutions $(\xi,\zeta)=(a,b)$ and $(\xi,\zeta)=(k-a,-b)$. This completes the proof of Statement~(D) as well as completes the proof of One Radius Theorem~\ref{Basic-theorem}, p~\pageref{Basic-theorem}.
\end{proof}
\end{proof}

\begin{proof}[\textbf{Proof of Corollary \ref{imaginary-part-vanishing-corollary}, p.~\pageref{imaginary-part-vanishing-corollary}.}]

According to Statement (A),
\begin{equation}\label{NYU-1}
    \int\limits_{S^k}\omega^{k/2+ib}dS_y=
    \int\limits_{S^k}\omega^{k/2-ib}dS_y\,.
\end{equation}
Comparing the imaginary parts of this identity  gives
\begin{equation}\label{imaginary-part-only}
\int\limits_{S^k}\omega^{k/2}\sin(b\ln\omega)dS_y=
    -\int\limits_{S^k}\omega^{k/2}\sin(b\ln\omega)dS_y\,,
\end{equation}
which implies that imaginary part is identically zero, and then, the proof of \eqref{imaginary-part-vanishing} is complete. The repeated differentiation with respect to $b$ of the imaginary part \eqref{imaginary-part-only} leads directly to \eqref{logarithm-integral}.

To proof the Taylor decomposition presented in \eqref{Taylor-Decomposition}, notice that $F(\alpha)$ is entire function and then for every $\alpha\in\mathbb{C}$,
\begin{equation}\label{Taylor-Deco-Proof-1}
    F(\alpha)=F(k/2)+F'(k/2)(\alpha-k/2)+\frac{F''(k/2)}{2!}(\alpha-k/2)^2+\cdots\,.
\end{equation}
The direct computation together with the first formula from \eqref{logarithm-integral}, p.~\pageref{logarithm-integral} yields
\begin{equation}\label{Taylor-Deco-Proof-2}
    F^{(2m+1)}(k/2)=\int\limits_{S^k}\omega^{k/2}(\ln\omega)^{2m+1}dS\equiv 0
\end{equation}
for $m=0,1,2,\ldots$ and
\begin{equation}\label{Taylor-Deco-Proof-3}
    F^{(2m)}(k/2)=\int\limits_{S^k}\omega^{k/2}(\ln\omega)^{2m}dS>0
\end{equation}
for $m=0,1,2,\cdots$. This completes the proof of the Corollary \ref{imaginary-part-vanishing-corollary}.
\end{proof}

\begin{proof}[Proof of Corollary \ref{One-dimension-sphe-pro} from p.~\pageref{One-dimension-sphe-pro}.]

Note that the condition: $a>b>0$ is sufficient to find $R$ and $r$
such that $R>r>0$ and the equalities $R^2+r^2=a\, ,\,\, 2R\, r=b$
hold. Note also that if we integrate a $2\pi$-periodic function
over the period from 0 to $2\pi$, we can replace $\sin(\theta)$ by
$\cos(\theta)$ and vice versa without changing the result of the
integration. Thus, the function under the integral can be reduced to the
Poisson kernel, which yields

\begin{equation}\label{22}
\begin{split}
& \int\limits_{0}^{2\pi} \frac{d\theta}{(a-b\,\sin(\theta))^p}
=\int\limits_{0}^{2\pi} \frac{d\theta}{(a+b\,\cos(\theta))^p}
\\& =(R^2-r^2)^{-p}\,\int\limits_{0}^{2\pi}
\left(\frac{R^2-r^2}{R^2+r^2+2R\,r\,\cos(\theta)}\right)^p\,d\theta
\\& =(a^2-b^2)^{-p/2}\,\int\limits_0^{2\pi}
\omega^p\, d\theta=
(a^2-b^2)^{-p/2}\int\limits_0^{2\pi}
\omega^{1-p}\, d\theta
\\& =(a^2-b^2)^{1/2-p}\, \int\limits_0^{2\pi}
(a-b\, \sin(\theta))^{p-1}\, d\theta\,,
\end{split}
\end{equation}
where the fourth equality was obtained by \eqref{Sphe-Pro-Statement-2}. This completes the proof of Corollary \ref{One-dimension-sphe-pro}.
\end{proof}

\begin{proof}[Proof of Corollary \ref{Basic-Th-Integrable-Case-general-Item} from p.~\pageref{Basic-Th-Integrable-Case-general-Item}.]

We organize the proof as the chain of identities, which will be explained at each step.

\begin{description}
  \item[Step 1.] Using the Statement (A) from p.~\pageref{Sphe-pro-direct-Statement} and comparing imaginary parts for both integrals with $\alpha=1+ib$ and $\beta=1-ib$, we have
    \begin{equation}\label{Step-1}
        \int\limits_{S^2(R)}\omega^{1+ib}dS_y=\int\limits_{S^2(R)}\omega\cos(b\ln\omega)dS_y.
    \end{equation}
  \item[Step 2.] Let $\Sigma$ be the unit sphere centered at point $x$. Then, we can observe that
     \begin{equation}\label{Step-2-1}
        dS_y=\frac{2R}{l+q}q^2d\Sigma_{\widetilde{y}}\quad\text{and}\quad
        d\psi=\frac{l+q}{4r\sin\psi}d\ln\omega\,,
    \end{equation}
    where all notations are pictured on the Figure \ref{Reference-Picture} below that can be used for a reference. Thus, the last integral in \eqref{Step-1} can be written as
    \begin{equation}\label{Step-2-2}
        \int\limits_{S^k(R)}\omega\cos(b\ln\omega)\frac{2R}{l+q}q^2d\Sigma_{\widetilde{y}}.
    \end{equation}

    \begin{figure}[!h]
    \centering
    \epsfig{figure=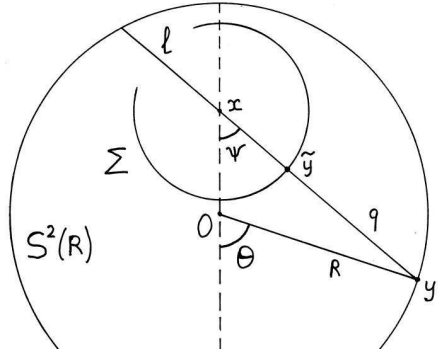,height=4.8cm}\\
  \caption{Reference picture.}\label{Reference-Picture}
    \end{figure}

  \item[Step 3.] Recall that according to formula~\eqref{omega-definition} from p.~\pageref{omega-definition} and Pythagorean theorem, we have
    \begin{equation}\label{Step-3-1}
        \omega(x,y)=\frac lq\quad\text{and}\quad lq=|R^2-r^2|.
    \end{equation}
    Note also that the integrand in \eqref{Step-2-2} depends only on angle $\psi$ pictured on Figure \ref{Reference-Picture} above. Therefore, the integral in \eqref{Step-2-2} becomes
    \begin{equation}\label{Step-3-2}
        4\pi R|R^2-r^2|\int\limits_{0}^{\pi}\frac{\sin\psi}{l+q}\cos(b\ln\omega)d\psi.
    \end{equation}
   \item[Step 4.] Observe that the last integrand is symmetric with respect to $\psi=\pi/2$. Indeed, $\sin(\pi/2+\tau)=\sin(\pi/2-\tau)$ and
\begin{equation}\label{Three_D_Dirichlet-29}
\begin{split}
    & (l+q)\circ\left(\frac{\pi}{2}-\tau\right)=(l+q)\circ\left(\frac{\pi}{2}+\tau\right)
    \\& \text{and}\quad\ln\frac{l}{q}\circ\left(\frac{\pi}{2}-\tau\right)
    =-\ln\frac{l}{q}\circ\left(\frac{\pi}{2}+\tau\right)\,.
\end{split}
\end{equation}
where the symbol $\circ$ denotes composition. Hence, the integral in \eqref{Step-3-2} can be written as
    \begin{equation}\label{Step-4-2}
        8\pi R|R^2-r^2|\int\limits_{0}^{\pi/2}\frac{\sin\psi}{l+q}\cos\left(b\ln\frac{l}{q}\right)d\psi.
    \end{equation}





   \item[Step 5.] The substitution for $d\psi$ presented in \eqref{Step-2-1},
   allows to rewrite the integral in \eqref{Step-4-2} as follows
    \begin{equation}\label{Step-5-2}
        \frac{8\pi R|R^2-r^2|}{4r}
        \int\limits_{\psi=0}^{\psi=\pi/2}\cos(b\ln\omega)d(\ln\omega).
    \end{equation}
   \item[Step 6.] The last integral can be computed directly and then the integral in \eqref{Step-5-2} yields
    \begin{equation}\label{Step-6-1}
        \left.\frac{2\pi R|R^2-r^2|}{br}\sin(b\ln\omega)\right|_{\psi=0}^{\psi=\pi/2}.
    \end{equation}
   \item[Step 7.] Note now that
    \begin{equation}\label{Step-7-1}
        \omega(\pi/2)=\frac{l(\pi/2)}{q(\pi/2)}=1\,,\quad\text{while}\quad\omega(0)=\frac{|R-r|}{R+r}.
    \end{equation}
    Therefore, the direct computation in \eqref{Step-6-1} gives
    \begin{equation}\label{Step-7-2}
        \frac{2\pi R|R^2-r^2|}{br}\sin\left(b\ln\frac{R+r}{|R-r|}\right)\,,
    \end{equation}
    which completes the proof of \eqref{Basic-Th-Integrable-Case-1}, p.~\pageref{Basic-Th-Integrable-Case-1} of the Corollary.
\end{description}
The partial case stated in \eqref{Basic-Th-Integrable-Case-2}, p.~\pageref{Basic-Th-Integrable-Case-2} can be obtained by taking limit at both parts of \eqref{Basic-Th-Integrable-Case-1} as $b\rightarrow0$. This limit exists since the integrand in \eqref{Basic-Th-Integrable-Case-1}, p.~\pageref{Basic-Th-Integrable-Case-1} or equivalently, in \eqref{Step-1}, p.~\pageref{Step-1} converges uniformly as $b\rightarrow0$. This completes the proof of Corollary~\ref{Basic-Th-Integrable-Case-general-Item}, p.~\pageref{Basic-Th-Integrable-Case-general-Item}.
\end{proof}


\section{Appendix}

The main goal of this section is to use the ratio of two segments $l$ and $q$ defined above, to analyze the behavior of $\omega(x,y)$ as $x$ approaches $y$. We show that $\lim\limits_{x\rightarrow y}\omega(x,y)\in[0,\infty]$ and depends on the path chosen for $x$ to approach $y$.

\begin{lemma}\label{omega_constant-lemma} Let $T$ be any point on the line $Oy$ such that $T\neq y$ and let $S^k_T(|Ty|)$ be the sphere of radius $\delta=|Ty|$ centered at $T$, where $|Ty|=|T-y|$ is the distance from $T$ to $y$. Then
\begin{equation}\label{Omega-constanta}
    \omega(x,y)=\frac{|OT|}{\delta}\quad\text{for all}\,\,\,x\in S^k_T(|Ty|)\setminus \{y\}\,.
\end{equation}
\end{lemma}

\begin{proof}[Proof of Lemma \ref{omega_constant-lemma}]
Let us choose $T\in\text{line}(Oy)$ as it is shown on Figure~\ref{Omega-constant} below.

\begin{figure}[!h]
    \centering
    \epsfig{figure=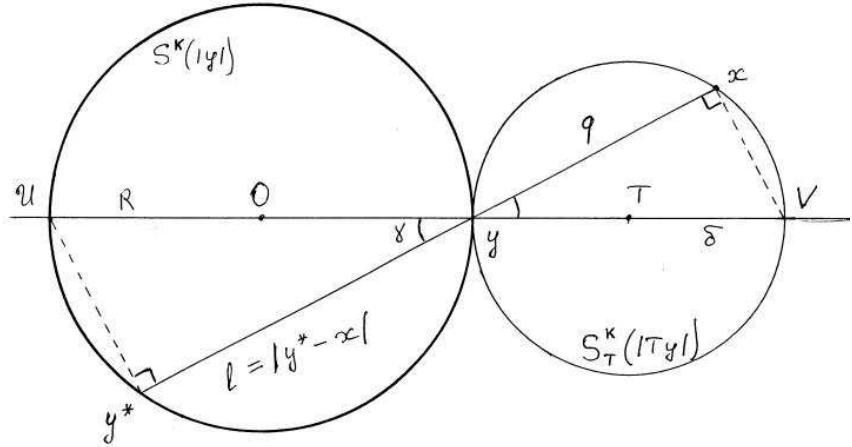,height=6cm}
    \caption{Level hypersurface for $\omega$.}\label{Omega-constant}
\end{figure}

Let $|Oy|=R$, $U$ and $V$ be the intersections of line $Oy$ with the spheres $S^k(|y|)$ and $S^k_T(|Ty|)$ respectively. $y^*\in S^k(|y|)$ is defined by the definition~\ref{Spherical_Ratio_Definition}. Let $\gamma=\angle Oyy^*=\angle xyT$. Observe that $\angle Uy^*y=\angle yxV=\pi/2$ since the segments $Uy$ and $yV$ are diameters. Therefore,
\begin{equation}\label{l-for-level-hypersurface}
    l=|y^*-x|=|y^*y|+|yx|=2R\cos\gamma+2\delta\cos\gamma
\end{equation}
and
\begin{equation}\label{q-for-level-hypersurface}
    q=|y-x|=2\delta\cos\gamma\,.
\end{equation}
Hence, the combination of \eqref{omega-definition}, \eqref{l-for-level-hypersurface} and \eqref{q-for-level-hypersurface} yields
\begin{equation}
    \omega(x,y)=\frac{l}{q}=1+\frac{R}{\delta}=\frac{\delta+R}{\delta}=\frac{|OT|}{\delta}\,.
\end{equation}
Therefore, for the chosen point $T$ formula \eqref{Omega-constanta} is justified. A similar argument yields \eqref{Omega-constanta} for all other positions of $T$ on the line $Oy$. This completes the proof of Lemma~\ref{omega_constant-lemma}.
\end{proof}

\begin{corollary}\label{limit-for-omega}
    Let $p$ be the $\mathcal{C}^1$ path through $y$ and let $x\in p$. Then
\begin{equation}\label{Omega-limit-formula}
    \lim\limits_{p\ni x\rightarrow y} \omega(x,y)= A(p)\in[0, \infty]\,,
\end{equation}
where the value $A(p)$ depends on the path $p$ chosen and can assume any value from the closed interval $[0,\infty]$.
\end{corollary}

\begin{proof}[Proof of Corollary \ref{limit-for-omega}]
Note first, if $p\subseteq S^k(|y|)$, then $\omega(x,y)\equiv0$ for every $x\in p\setminus \{y\}$ and therefore, $A(p)=0$. If $p\subseteq S^k_T(|yT|)$, then, by \eqref{Omega-constanta},
\begin{equation}
    A(p)=\frac{|OT|}{|\delta|}=\frac{|OT|}{|yT|}\,,
\end{equation}
which can be any number from $(0,\infty)$ depending on $T$ chosen on the line $Oy$. In particular, if $T=\infty$, then $A(p)=1$, which corresponds to the case when $p$ belongs to the tangent hyperplane to $S^k(|y|)$ at point $y$. In this case $\omega(x,y)\equiv1$ for every $x\in p\setminus \{y\}$. Finally, to get $A(p)=\infty$, $p$ must be a $\mathcal{C}^1$ curve non-tangential to $S^k(|y|)$ at point $y$. Indeed, let $\gamma_0$ be the angle between $p$ and line$(Oy)$ at point $y$ as it is pictured on Figure~\ref{Omega-Limit-picture} below.
\begin{figure}[!h]
    \centering
    \epsfig{figure=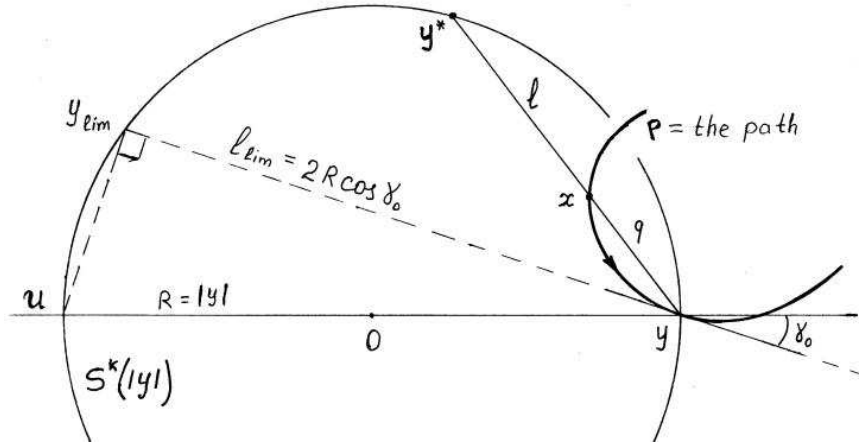,height=6cm}
    \caption{Non-tangential limit for $\omega$.}\label{Omega-Limit-picture}
\end{figure}

Note that $Uy$ is the diameter which implies that $\angle Uy_{\lim}y=\pi/2$, and then,
\begin{equation}
    \lim\limits_{p\ni x\rightarrow y} l(x)=l_{\lim}=2R\cos\gamma_0\,,
\end{equation}
while
\begin{equation}
    \lim\limits_{p\ni x\rightarrow y} q(x)=\lim\limits_{p\ni x\rightarrow y} |x-y|=0\,.
\end{equation}
Therefore,
\begin{equation}
    \lim\limits_{p\ni x\rightarrow y} \omega(x,y)= \frac{l(x)}{q(x)}=\infty\,,
\end{equation}
which completes the proof of Corollary \eqref{limit-for-omega}.
\end{proof}

\bibliographystyle{IEEEtran}

\end{document}